\documentclass[12pt,a4paper]{article}
\usepackage[T1]{fontenc}
\usepackage[latin1]{inputenc}
\usepackage{amssymb}
\usepackage{latexsym}
\usepackage{amsmath, amsthm, url, graphicx}

\usepackage{graphicx}

\newcommand{\bbQ}[0]{{\mathbb Q}}
\newcommand{\bbR}[0]{{\mathbb R}}

\newcommand{\posrat}[0]{{\mathbb Q}_+}

\newcommand{\ve}[0]{{\varepsilon}}

\newcommand{\cov}[0]{{\, \lhd \,}}
\newcommand{\kov}[0]{{\, \lessdot \,}}

\newcommand{\balcov}[0]{\sqsubseteq}
\newcommand{\bal}[0]{{\sf b}}

\newcommand{\formal}[1]{{\cal #1}}
\newcommand{\radius}[1]{\rho(#1)}
\newcommand{\complete}[1]{\ensuremath{\formal{M}(#1)}}
\newcommand{\pos}{\ensuremath{\operatorname{Pos}}}
\newcommand{\Pos}[1]{\ensuremath{\pos( #1)}}
\newcommand{\compr}[2]{\ensuremath{\{#1  \mid  #2 \}}}
\newcommand{\lowp}[1]{\ensuremath{\operatorname{Pos} #1}}
\newcommand{\powfe}{{\cal P}_{\text{K-fin}}}

\newcommand{\intersects}[2]{#1 \between #2}

\newcommand{\Pt}{\ensuremath{\operatorname{Pt}}}

\newtheorem{thm}{Theorem}[section]
\newtheorem{cor}[thm]{Corollary}
\newtheorem{prop}[thm]{Proposition}
\newtheorem{lemma}[thm]{Lemma}

\newtheorem{defin}[thm]{Definition}
\newtheorem{example}[thm]{Example}

\newcommand{\AMSclass}[2]{{\textbf{#1 A.M.S. subject classification:} #2}}
\newcommand{\keywords}[1]{{\textbf{Keywords:} #1}}

\title{Metric complements of overt closed sets}
\author{Thierry Coquand%
\footnote{Computer science and engineering,
Chalmers University of Technology and University of Gothenburg,
SE-412 96 G\"oteborg, Sweden.}
\and Erik Palmgren\footnote{Swedish Collegium for Advanced Study
S-752 38 Uppsala, Sweden and\newline
Department of Mathematics,
Uppsala University,
P.O. Box 480,
S-751 06 Uppsala, Sweden.\newline
Partly supported by a grant from the Swedish Research Council (VR).}
\and Bas Spitters%
\footnote{Eindhoven University of Technology, Department of
Mathematics and Computer Science, P.O. Box 513, 5600 MB Eindhoven, the
Netherlands. Supported by NWO.}}
\date{}
\begin{document}
\maketitle
\abstract{We show that the set of points of an overt closed subspace of a
metric completion of a Bishop-locally compact metric space is located.
Consequently, if the subspace is, moreover,
compact, then its collection of points is Bishop compact.}


\AMSclass{2000}{06D22, 
28C05}

\section{Introduction}
We continue our investigations into the relations between Bishop's constructive
mathematics and formal
topology~\cite{Coquand:jucs_11_12:formal_topotoly_and_constructive,
Palmgren:crealform,Palmgren:metric,Located-overt,integrals-valuations,
Palmgren:OpenSub}.
Previously, we gave a formal definition of locatedness~\cite{Located-overt} and
showed that an overt closed subspace of a compact
formal space is (formally) located. Here we consider a generalization of the
pointwise side of this result. We use the real numbers as a running example.
We work in informal Bishop-style mathematics, including the axiom of dependent
choice.

\section{Preliminaries}
\subsection*{Bishop}
We assume familarity with Bishop's constructive
mathematics~\cite{Bishop/Bridges:1985}, but we
recall some relevant notions.

\begin{defin}
  A set is \emph{finite} if it is in bijective correspondence with a set
  $\{0, \ldots, n\}$, $n \geqslant 0$. A set is {\emph{Kuratowski finite}}
  ({\emph{K-finite}}, \emph{finitely enumerable}) if it is the image of a
finite set.
\end{defin}
A subset of a set $X$ if K-finite iff it can written as
$\{x_1,\ldots,x_n\}$ with $x_1,\ldots,x_n$ in $X$. Such a set does not need to
have a cardinality if the equality on $X$ is not decidable. For example, the
set $\{a,b\}$ is K-finite. However, it is finite iff we can decide whether
$a=b$.

\begin{defin}
A metric space is said to be {\emph{totally bounded}} if for each $\varepsilon
> 0$ the space can be covered by a K-finite set of balls with
radius at most $\varepsilon$. A subset of a metric space is
\emph{Bishop-compact} if it is complete and totally bounded.\\
A metric space is said to be {\emph{locally totally
bounded}} if for each ball and each $\varepsilon > 0$ the ball can be covered
by a K-finite set of balls with radius at most $\varepsilon$. A metric space is
\emph{Bishop-locally compact} if it is complete and locally totally
bounded.
\end{defin}

The closed unit interval is compact. The real numbers are locally compact.

\begin{defin}
  \label{def:located-metric}A subset $A$ of a metric space $(X, \rho)$ is
  located if for each $x$ in $X$ the distance $\inf \{\rho (x, a)\mid a \in
  A\}$ exists as a (Dedekind) real number.
\end{defin}

In classical mathematics all sets are located. Constructively this is not the
case, as the following Brouwerian counterexample shows.
\begin{example}\label{ex:not-located}
Consider the set
\[\compr{x \in \bbR}{x > 1 \text{ or } (x > 0 \text{ and } P)}.\]
This set will only be located if we can decide whether the proposition $P$
holds.
\end{example}

\begin{defin}
A subset of a metric space is {\emph{Bishop-closed}} if it contains all its
limit points, i.e.~if it coincides with its closure.
\end{defin}

The closed unit interval [0,1] is Bishop-closed.

A Bishop closed located
subset of a metric space coincides with the complement of its complement: a
Bishop closed located set coincides with the set of all points which have zero
distance to it.

\subsection*{Formal topology}
\begin{defin}
A {\emph{formal topology}}~{\cite{Sambin:SomePoints}}
consists of a pre-order $(S, \leqslant)$ of basic opens and a
relation $\vartriangleleft \subset S \times \mathcal{P}(S)$, the covering
relation, which satisfies:
\begin{description}
  \item[Ref] $a \in U$ implies $a \vartriangleleft U$;
  \item[Tra] $a \vartriangleleft U,$ $U \vartriangleleft V$ implies $a
  \vartriangleleft V$, where $U \vartriangleleft V$ means $u \vartriangleleft
  V$ for all $u \in U$;
  \item[Loc] $a \vartriangleleft U$, $a \vartriangleleft V$ implies $a
  \vartriangleleft U \wedge V = \{x\mid \exists u \in U \exists v \in V.x
  \leqslant u, x \leqslant v\}$;
  \item[Ext] $a \leqslant b$ implies $a \vartriangleleft \{b\}$.
\end{description}
These axioms are known as Reflexivity, Transitivity, Localization and
Extensionality. {\textbf{Ref}} and {\textbf{Ext}} say that if a basic open
belongs to a family, then the family covers it. {\textbf{Tra}} is the
transitivity of the cover. {\textbf{Loc}} is the distributive rule for
frames.

The {\emph{formal intersection}} $U \wedge V$ is defined as $U_{\leqslant}
\cap V_{\leqslant}$, where $Z_{\leqslant}$ is the set $\{x\mid \exists z \in Z.x
\leqslant z\}$. Another common notation for $Z_{\leqslant}$ is
$Z_{\downarrow}$. We write $a \vartriangleleft b$ for $a \vartriangleleft
\{b\}$. We write $U \equiv V$ iff $U \cov V$ and $V \cov U$.
\end{defin}

\begin{defin}
  Let $(S, \cov)$ be a formal topology. A \emph{point} is an inhabited subset
$\alpha
  \subset S$ which is filtering with respect to $\leqslant$, and such that $U
  \cap \alpha$ is inhabited, whenever $a \cov U$ for some $a \in \alpha$. The
  collection of points is denoted by $\Pt (S)$. Let $U$ be an open in $S$. Then
$U_*$ denotes the class of points $\alpha$ such that $\alpha\in U$.
\end{defin}

\begin{example}
  \label{ex:locale-reals}The formal reals are inductively defined by the
following relation on the open rational intervals ordered by inclusion.
  \begin{enumerate}
    \item $(p, s) \cov \{(p, r), (q, s)\}$ if $p \leqslant q < r \leqslant s$;

    \item $(p, q) \cov \compr{(p', q')}{p < p' < q' < q}$.
  \end{enumerate}
  The points of this space are precisely the (Dedekind) real numbers.
\end{example}

\begin{defin}
Let $(S, \vartriangleleft)$ be a
formal topology. A \emph{sublocale} is a formal topology $(S, \cov')$ such
that $\vartriangleleft \subset \vartriangleleft'$ and $a \wedge' b
\vartriangleleft' a \wedge b$.\\
Let $U \subset S$. The \emph{closed}
sublocale $S \setminus U$ is $u \vartriangleleft_{- U} V$ iff $u
\vartriangleleft V \cup U$.
\end{defin}

\begin{example}
  \label{ex:closed-interval}The set $\compr{(p, q)}{q \leqslant 0 \vee p
  \geqslant 1}$ represents the closed unit interval as a subspace of the real
  line.
\end{example}

\begin{defin}
  $\pos$ is called a {\emph{positivity predicate}} on a formal topology $S$ if
  it satisfies:
  \begin{description}
    \item[Pos] $U \vartriangleleft U^+$, where $U^+ := \{u \in U\mid \Pos{u}
    \}$.
    \item[Mon] If $\Pos{u}$ and $u \vartriangleleft V$, then $\Pos{V}$
    --- that is, $\Pos{v}$ for some $v \in V$.
  \end{description}
 A formal space is \emph{overt} if it carries a positivity predicate.
\end{defin}
Impredicatively, a formal topology is overt iff the locale it generates is
overt, or open.

Classically, all formal topologies are overt. Constructively this is not the
case, as the following formal analogue of Example~\ref{ex:not-located} shows.
\begin{example}
  \label{ex:not-overt} The closed sublocale defined by the open
  \[ \compr{(p, 0)}{p < 0} \cup \compr{(2, q)}{q > 2} \cup \compr{(1, q)}{q >
     1 \text{ and } P} \]
  is overt if we can decide whether the proposition $P$ holds;
see~\cite{Located-overt}.
\end{example}

\subsection*{Metric completion}
\begin{defin}
  To any metric space $X$, we define, following Vickers~{\cite{Vickers:locAA}}
  and Palmgren~{\cite{Palmgren:metric}}, a formal topology $\formal{M}(X)$
called the {\emph{localic completion}} of $X$. A formal open is a pair $(x, r)
\in X
  \times \bbQ^{> 0}$, written $\bal(x,r)$. We define the relation $\bal(x,r) <
\bal(y,s)$ iff $d (x, y) < s - r$ as illustrated below.
\begin{center}
    \includegraphics[width= 5cm]{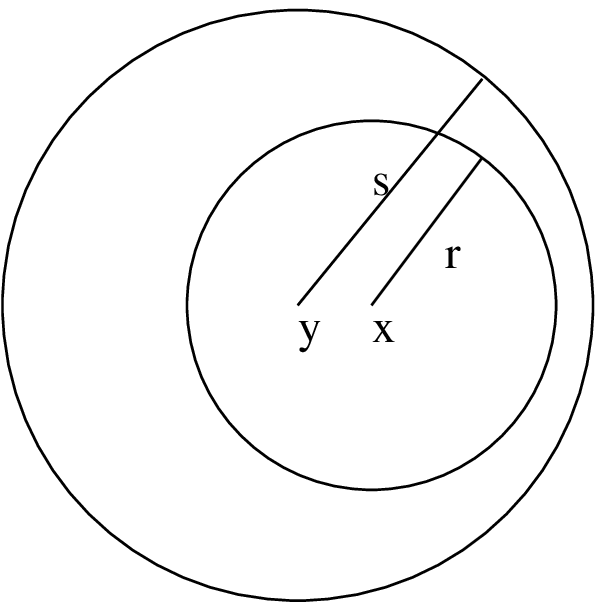}
\end{center}

  The order $\leqslant$ is defined by $\bal(x,r) \leqslant \bal(y,s)$ iff
$d (x,  y) < t$ for all $t > s - r$. The covering relation $\cov$ is
inductively generated by the axioms

  \begin{description}
    \item[M1] $u \vartriangleleft \{v\mid v < u\}$;
    \item[M2] $\complete{X} \vartriangleleft \{\bal(x,r)\mid x \in X\}$
for any $r$.
  \end{description}

  {\textbf{M1}}: Every ball is covered by all the balls strictly inside it
  (since the ball is open). {\textbf{M2}}: For each $r > 0$, the space is
  covered by all balls of size $r$.\\
  We define $U<V:=\forall u\in U\exists v\in V.u<v$.
\end{defin}

\begin{example}
  Consider the formal unit interval $[0, 1]$. Then $[0, 1] = \bal(0,3) =
\bal(0,2)$, but $\bal(0,3) > \bal(0,2)$.

  Similarly, $\bal(0,3) \cov\bal(0,2)$, but it is not the case
that $\bal(0,3) \leqslant\bal(0,2)$. This shows that $a \cov b$ does not imply
$a\leqslant b$.
\end{example}

\begin{prop}
  The localic completion of a metric space is always overt.
\end{prop}

The formal reals are the metric completion of the rational numbers.

\subsection*{Elementary description of the cover of  $\formal{M}(X)$}
The $\kov$ cover relation below, introduced in~\cite{Palmgren:metric},
generalises
the one introduced by Vermeulen~\cite{Vermeulen} and
Coquand~\cite{CederquistNegri} for ${\mathbb R}$:
\begin{eqnarray*}
p \balcov_{\ve} U &:=& (\forall q \le p)[{\rm radius}(q)
\le \ve \Rightarrow
\{q\} \le U]\\
p \balcov U&:=& p \balcov_{\ve} U\text{ for some }\ve \in \posrat\\
A(b,c) &:=&\{C \in \powfe(\formal{M}_X) : b \balcov C < c\}\\
a\kov U &:=& (\forall b < c < a)\,(\exists U_0 \in A(b,c))\, U_0 < U.
\end{eqnarray*}

\begin{thm}[\cite{Palmgren:metric}]
If $X$ is a Bishop locally compact metric space, then
\[a \cov U \Longleftrightarrow a \kov U\]
($\Leftarrow$ holds for any metric $X$ space.)
\end{thm}

\begin{thm}[\cite{Palmgren:metric}] Let $X$ be a complete metric
space. Then there exist a metric isomorphism $j$ between $X$ and
$\Pt(\formal{M}(X))$.
\end{thm}

In particular, this holds for the real numbers.

\section{Main results}
We write $B(x,r)$ for the set $\{y\mid d(x,y)<r\}$. Then $\bal(x,r)_*=B(x,r)$.
\begin{lemma}\label{lm1}
Let $X$ be a metric space. Then an inhabited set $S \subseteq X$
is located if, and only if, for all $x \in X$ and all positive $\delta <
\ve$ we have
\[
S \cap B(x,\delta)= \emptyset\text{ or }\intersects{S}{B(x,\ve)}.
\]
\end{lemma}
Where $\intersects{A}{B}$ means that $A \cap B$ is inhabited.
\begin{lemma}\label{lm2}
Let $X$ be a metric space and let $M= \formal{M}(X)$ be its localic completion.
If $O
\subseteq M$ and the sublocale $M \setminus O$ is overt, then any positive
neighbourhood contains a point of $M\setminus O$.
\end{lemma}
\begin{proof}
Suppose that $P$ is the positivity predicate of $M
\setminus O$.
Denote the cover relation of $M \setminus O$ by $\cov'$.

Suppose $a=\bal(x,\delta) \in P$. Let $a_1=a$.
Suppose we have constructed in $P$:
$$a_1 \ge a_2 \ge \cdots \ge a_n,$$
so that ${\rm radius}(a_{k+1}) \le {\rm radius}(a_k)/2$.

By (M1) and localisation we get
$$a_n \cov' \{a_n\} \wedge \{\bal(y,\rho) : y \in X\}$$
where $\rho = {\rm radius}(a_n)/2$.
Since $a_n \in P$ we obtain some $b \in \{a_n\} \wedge \{\bal(y,\rho) : y \in
X\}$
with
$b \in P$. Clearly ${\rm radius}(b) \le  {\rm radius}(a_n)/2$. Let $a_{n+1} =b$.

Let
$$\alpha = \{ p\in M : (\exists n) a_n \le p \}.$$
Since the radii of $a_n$ are shrinking, this defines a point in ${\rm Pt}(M)$.
 (Note that we used Dependent Choice.)

We claim that $\alpha \in  {\rm Pt}(M \setminus O) = {\rm Pt}(M) \setminus O_*$.
Suppose that $\alpha \in O_*$, i.e.\ for some $c \in O$: $c \in \alpha$. Hence
there
is $n$ with $a_n \le c$.  Thus $a_n \cov O$, that is $a_n \cov' \emptyset$. But
since $a_n$ is positive, this is impossible! So $\alpha \in  {\rm Pt}(M
\setminus O)$.
\end{proof}

The following theorem can be conveniently formulated using the following
definition. However, no futher facts about this definition are needed.
\begin{defin}\cite{Located-overt}\label{def:located}
  Let $X$ be a metric space. A predicate
  $\lowp{}$ on $S = \compr{\bal(x,r)}{x \in X, r \in \bbQ^+}$ is called
  {\emph{located}} if
\begin{itemize}
  \item $\lowp (u)$ and $u \cov V$ imply that $\lowp (v)$ for some $v$ in $V$;
  \item $v < u$ implies that $\neg \lowp{v}$ or $\lowp{u}$.
\end{itemize}
  Let $T$ be a closed sublocale of $\complete{X}$. Then $T$ is called
  {\emph{located}} if there is a located predicate {\lowp{}} such that $T$
  coincides with the closed sublocale defined by the open $\neg \lowp{}
  \subset S$.
\end{defin}

\begin{thm}\label{thm1}
Let $X$ be a Bishop locally compact metric space and let $M= \formal{M}(X)$ be
its localic completion. Let $O \subseteq M$. If a sublocale $M \setminus O$ is
overt, then $M\setminus O$ is (formally) located. Consequently, the set of
points $Y= j^{-1}[{\rm Pt}(M \setminus O)]$ is located as a subset of $X$ and
moreover $O_*$ is the metric complement of $Y$ in $X$.
\end{thm}
\begin{proof}
That $M \setminus O$ is overt means that
 there is an inhabited  subset $P \subseteq M$ (a set of positive formal
neighbourhoods) so that
\begin{itemize}
\item[(P1)] $a \cov_M\; O \cup U$ and $a \in P$ implies
$U \cap P$ inhabited,
\item[(P2)] $U \cov_M\; O \cup (U \cap P)$.
\end{itemize}

Now since $P$ is inhabited, Lemma~\ref{lm2} ensures that $Y$
is inhabited. Let $x \in X$ be arbitrary. Consider positive rational
numbers $\delta < \ve$. Take $\ve'$ with $\delta < \ve' < \ve$
and let $\theta = \ve - \ve'$. Using the (P2), (M2) and the
localisation we get
$$\bal(x,\ve) \cov_M O \cup (P \cap \{c \in M : \radius{c} = \theta/2\}).$$
Thus also
$$\bal(x,\ve) \kov O \cup (P \cap \{c \in M : \radius{c} = \theta/2\}).$$
and by definition there is a K-finite $W \in A(\bal(x,\delta),\bal(x,\ve'))$
with
$$W < O \cup (P \cap \{c \in M : \radius{c} = \theta/2\}).$$
Since $W$ is K-finite we have one of the cases
\begin{itemize}
\item[(C1)] $W< O$,
\item[(C2)] $(\exists d \in W) d < P \cap  \{c \in M : \radius{c} = \theta/2\}$.
\end{itemize}

In case (C1) we have $\bal(x,\delta) \cov O$ and hence
 $\bal(x,\delta)_* \subseteq O_*$. Thus $\bal(x,\delta)_*\, \cap Y = \emptyset$.

 In case (C2) there is $d \in W$ and $c \in P$ with  $d < c$ and $\radius{c}=\theta/2$. Suppose $c= \bal(y,\theta/2)$ and
 $d= \bal(z,\tau)$. Now $W < \bal(x,\ve')$. Hence $d(z,x) < \ve'$.
Moreover $d < c$ implies $d(z,y) +\tau <\theta/2$. Thus
\[d(x,y) \le d(x,z) + d(z,y) < \ve' + \theta/2-\tau = \ve -\theta/2 -\tau < \ve
-\theta/2.\]
Thereby $c < \bal(x,\ve)$, and so $\bal(x,\ve)\in P$. By Lemma~\ref{lm1},
$\intersects{Y}{B(x,\ve)}$.

We have thus showed that $M\setminus O$, and hence $Y$, is located. Using
Lemma~\ref{lm2}, we have that
\[d(x,Y)>0 \Longleftrightarrow (\exists \delta >0) B(x,\delta) \cap Y
=\emptyset.\]
We claim that $O_*$ is the metric complement of $Y$, i.e.
$$ x \in O_* \Longleftrightarrow d(x,Y) >0.$$
If $x \in O_*$ then for some $\delta>0$, $B(x,\delta) \subseteq O_*$. Thus
$B(x,\delta) \cap Y$
cannot be inhabited. Conversely, suppose that $B(x,\delta) \cap Y =\emptyset$
for some $\delta>0$.
We have by (P2), that
$$\bal(x,\delta) \cov O \cup (\{\bal(x,\delta)\} \cap P)$$
Thus $B(x,\delta) \subseteq O_* \cup (\{\bal(x,\delta)\} \cap P)_*$, and hence $x \in O_*$ or
$x \in  (\{\bal(x,\delta)\} \cap P)_*$. In the latter case $\bal(x,\delta) \in P$, which contradicts
$B(x,\delta) \cap Y =\emptyset$. Thus $x \in O_*$.
\end{proof}

\begin{thm}
Let $X$ be a metric space and let $M= \formal{M}(X)$ be its localic completion.
Let $O \subseteq M$ be such that $M\setminus O$ is compact and overt.
Then ${\rm Pt}(M \setminus O)$ is Bishop-compact.
\end{thm}
\begin{proof}
Let $P$ and $\cov'$ be as in the proof of Lemma~\ref{lm2}.
Let $\ve>0$ be given. Then by axiom M2 and positivity
$$M \cov' \{\bal(x,\ve/2) : x \in X \} \cov'  \{\bal(x,\ve/2) : x \in X \}  \cap
P.$$
By compactness there is some K-finite
\[F =\{\bal(x_1,\ve/2),\ldots,
\bal(x_n,\ve/2)\}
\subseteq  \{\bal(x,\ve/2) : x \in X \}  \cap P\]
so that
\begin{equation}\label{poscov}
M \cov' F.
\end{equation}
Since each $\bal(x_i,\ve/2)$ is positive there is by Lemma~\ref{lm2} some
$\alpha_i \in
\bal(x_i,\ve/2)_*$
which is in ${\rm Pt}(M\setminus O)$. By (\ref{poscov}), each
point in ${\rm Pt}(M\setminus O)$
has distance smaller than $\ve$ to some point $\alpha_i$. Thus
$\{\alpha_1,\ldots,\alpha_n\}$
is the required $\ve$-net.
\end{proof}

\begin{cor}
If in the context of Theorem~\ref{thm1},  $X$ is Bishop-compact, and then so is
$Y$.
\end{cor}
\begin{proof}
 If $X$ is Bishop-compact, then $\formal{M}(X)$ is a compact as a formal
space~\cite{Vickers:locAA}. Hence, $M\setminus O$ is compact. By
Theorem~\ref{thm1}, $Y$ is Bishop-compact.
\end{proof}

\bibliographystyle{alpha}
\bibliography{compact,Erik,located}

\begin{thebibliography}{Sam03}

\bibitem[BB85]{Bishop/Bridges:1985}
Errett Bishop and Douglas Bridges.
\newblock {\em {Constructive analysis}}, volume 279 of {\em {Grundlehren der
  Mathematischen Wissenschaften}}.
\newblock Springer-Verlag, 1985.

\bibitem[CN96]{CederquistNegri}
Jan Cederquist and Sara Negri.
\newblock A constructive proof of the {H}eine-{B}orel covering theorem for
  formal reals.
\newblock In {\em Types for proofs and programs (Torino, 1995)}, volume 1158 of
  {\em Lecture Notes in Comput. Sci.}, pages 62--75. Springer, 1996.

\bibitem[CS05]{Coquand:jucs_11_12:formal_topotoly_and_constructive}
T.~Coquand and B.~Spitters.
\newblock Formal {T}opology and {C}onstructive {M}athematics: the {G}elfand and
  {S}tone-{Y}osida {R}epresentation {T}heorems.
\newblock {\em Journal of Universal Computer Science}, 11(12):1932--1944, 2005.

\bibitem[CS09]{integrals-valuations}
Thierry Coquand and Bas Spitters.
\newblock Integrals and valuations.
\newblock {\em Journal of Logic and Analysis}, 1(3):1--22, 2009.

\bibitem[Pal05]{Palmgren:crealform}
Erik Palmgren.
\newblock Continuity on the real line and in formal spaces.
\newblock In P.~Schuster L.~Crosilla, editor, {\em From {S}ets and {T}ypes to
  {T}opology and {A}nalysis: {T}owards {P}racticable {F}oundations for
  {C}onstructive {M}athematics}. Oxford UP, 2005.

\bibitem[Pal07]{Palmgren:metric}
Erik Palmgren.
\newblock A constructive and functorial embedding of locally compact metric
  spaces into locales.
\newblock {\em Topology Appl.}, 154(9):1854--1880, 2007.

\bibitem[Pal09]{Palmgren:OpenSub}
Erik Palmgren.
\newblock Open sublocales of localic completions.
\newblock Department of Mathematics Report 2009:1, 2009.

\bibitem[Sam03]{Sambin:SomePoints}
Giovanni Sambin.
\newblock Some points in formal topology.
\newblock {\em Theoret. Comput. Sci.}, 305(1-3):347--408, 2003.
\newblock Topology in computer science (Schlo\ss\ Dagstuhl, 2000).

\bibitem[Spi07]{Located-overt}
Bas Spitters.
\newblock Locatedness and overt sublocales.
\newblock Preprint available at \url{http://arxiv.org/abs/math/0703561}, 2007.

\bibitem[Ver86]{Vermeulen}
Japie Vermeulen.
\newblock Constructive techniques in functional analysis, 1986.

\bibitem[Vic05]{Vickers:locAA}
Steven Vickers.
\newblock Localic completion of generalized metric spaces. {I}.
\newblock {\em Theory Appl. Categ.}, 14:No. 15, 328--356, 2005.

\end{thebibliography}
\end{document}